\documentclass[a4paper,12pt]{article}
\usepackage[francais,english]{babel}
\usepackage{latexsym,amssymb}
\usepackage{amsmath}
\usepackage{pb-diagram}
\newtheorem{thm}{Theorem}[section]

\def\SO{{\rm{SO}}}
\def\Sp{{\rm{Sp}}}
\newcommand{\HH}{\mathbb{H}}
\newcommand{\C}{\mathbb{C}}
\newcommand{\Z}{\mathbb{Z}}
\newcommand{\R}{\mathbb{R}}
\newcommand{\PP}{\mathbb{P}}

\newcommand{\kker}{\mathrm{Ker} \,}
\newcommand{\lto}{\ensuremath{\longrightarrow}}
\newcommand{\function}[5]
{\begin{eqnarray*}\begin{array}{r@{}ccl}
#1\;\colon\;  & #2 &\lto & #3 \\[.05cm]
  & #4 &\longmapsto  & #5
\end{array}\end{eqnarray*}
}
\title{\bf{Eigenvalues of the Basic Dirac Operator on Quaternion-K\"ahler Foliations}}
\author{{\small {\bf {Georges Habib}}}\\
{\small Institut {\'E}lie Cartan,
Universit{\'e} Henri Poincar{\'e}, Nancy I, B.P. 239}\\
{\small 54506 Vand\oe uvre-L{\`e}s-Nancy Cedex, France}\\
{\small {\tt habib@iecn.u-nancy.fr}}}
\date{}
\begin{document}
\maketitle
\begin{abstract}
\noindent
In this paper, we give an optimal  lower bound for the eigenvalues of the basic Dirac operator  on a quaternion-K\"ahler foliations. The limiting case is characterized by the existence of  quaternion-K\"ahler Killing  spinors. We end this paper by giving some examples.
\end{abstract}
{\bf Key words:} Basic Dirac operator, quaternion-K\"ahler foliations,  eigenvalues, quaternion-K\"ahler Killing spinors.\\\\
{\bf Mathematics Subject Classification:}  53C20, 53C12, 57R30, 58G25
\section{Introduction}
On a compact quaternion-K\"ahler spin  manifold $(M,g)$ of dimension $4m\geq8,$ O. Hijazi and J.-L. Milhorat \cite{HM3} conjectured that any eigenvalue of the Dirac operator satisfies 
\begin{equation}
\lambda^2\geq\frac{m+3}{4(m+2)}S,
\label{eq:17}
\end{equation}
where $S$ denotes the constant scalar curvature (such manifolds are  Einstein \cite{Bes}). They proved that \eqref{eq:17} is   true for  $m=2$ and $m=3.$ For this, they introduced \cite{HM1} the {\it twistor operator}, as in the K\"ahler case,  on each eigenbundle associated with the eigenvalues of the fundamental $4$-form $\Omega$ \cite{HM2}.  Using  representation theory,  the lower bound  \eqref{eq:17} is established by W. Kramer, U. Semmelmann and G. Weingart \cite{ksw3}. Their proof is based on the decomposition in two ways of the bundle $TM\otimes TM \otimes \Sigma M$ into parallel subbundles under the action of the group $\Sp_{1} \times \Sp_{m}.$\\
On a compact Riemannian manifold $(M,g_{M},\mathcal{F})$ with a  spin foliation $\mathcal{F}$ of codimension $q$ and a bundle-like metric $g_{M},$  S. D. Jung \cite{J} gives a Friedrich-type inequality. For  K\"ahler foliations, he also gives a Kirchberg-type inequality for  odd complex dimensions  \cite{s} where the even case was proved by the author  \cite{GH}.  The main result of this paper is to prove the following theorem: 
\begin{thm}  Let $(M, g_{M},\mathcal{F})$ be a compact Riemannian manifold with a quaternion-K\"ahler spin foliation $\mathcal{F}$ of codimension $q=4m$ and  a bundle-like metric $g_{M}$  with a coclosed  basic $1$-form  mean curvature $\kappa.$ Then the foliation is minimal  and any eigenvalue $\lambda$ of the basic Dirac operator satisfies 
\begin{equation}
\lambda^{2}\geq\frac{m+3}{4(m+2)}\sigma^{\nabla}, 
\label{eq:16}
\end{equation}
where $\sigma^{\nabla}$ denotes the transversal scalar curvature.
\end{thm}
Our approach comes from an adaptation of  \cite{HM4} and  \cite{ksw2} to the case of Riemannian foliations where the key point is to prove that the mean curvature vanishes since the transversal Ricci curvature is strictly positive. The limiting case is characterized by the existence of quaternion-K\"ahler Killing spinors (see section \ref{sec:4} for  details).\\ 
  We point out that throughout this paper, we consider a bundle-like metric such that the mean curvature is a basic $1$-form and coclosed. The existence of such metric is assured in \cite{ Ma, MRM}. \\The author would like to thank J.-L. Milhorat for helpful discussions also he would like to thank Oussama Hijazi for his encouragment.
\section{ Spin Foliations}
In this section, we  summarize some standard facts about spin 
foliations. For details, we refer to  \cite{J}, \cite{gk1}, \cite{GH}, \cite{p}.\\\\
Let $(M,g_{M}, \mathcal{F})$ be a $(p+q)$-dimensional Riemannian  manifold
with a  Riemannian foliation  $\mathcal{F}$ of codimension $q$ and let  $\nabla^{M}$ be
the Levi-civita connection associated with $g_{M}.$  We denote by
 $L$  the tangent bundle of $TM$ and $Q= TM/L \simeq
L^\perp$ the normal bundle and we assume  $g_{M}$ to be a {\it bundle-like metric} on
$Q$, that means the induced metric $g_{Q}$ verifies for all $X\in  \Gamma(L)$ the holonomy invariance condition that is
 $ \mathcal{L}_{X}g_{Q}=0,$ where $\mathcal{L}_{X}$ is the Lie derivative with respect to $X$.
Let $\nabla $ be the transversal Levi-Civita connection on $Q$ defined for all $Y\in\Gamma(Q)$ by 
\begin{equation*}
\nabla _{X} Y =
\left\{\begin{array}{ll}
\pi [X,Y], &  \textrm {$ \qquad{\forall}  X  \in \Gamma(L)$ },\\
\pi (\nabla_{X}^{M}Y), & \textrm {$ \qquad{\forall} X  \in \Gamma (Q)$ },
\end{array}\right.
\end{equation*}
where $\pi: TM \rightarrow Q$ denotes the projection. The curvature of
 $\nabla $ acts on $\Gamma(Q)$ by :
$$R^{\nabla} \left(X,Y\right)=-\nabla _{X}\nabla _{Y}+\nabla _{Y}\nabla 
_{X}+\nabla _{\left[X,Y\right]},\qquad{\forall} X,Y \in \chi 
\left(M\right).$$
 We denote by $\rho^{\nabla}, \sigma^{\nabla} $   the transversal Ricci curvature  and  the scalar curvature  respectively associated with  $\nabla.$
The foliation  $\mathcal{F}$ is said to be transversally Einstein if and
only if  $\rho ^{\nabla }=\frac{1}{q}\sigma ^{\nabla } {\rm Id},$  with constant transversal scalar curvature.
The mean curvature of $\mathcal{F}$ is given for all $X \in 
\Gamma(Q)$  by  $\kappa \left(X\right)= g_{Q}\left(\tau ,X\right)$,  where $\tau $ is the trace of the second fundamental form $II$ of  $\mathcal{F}$ 
defined  by:
\function {II}{\Gamma(L) \times \Gamma(L) }{\Gamma(Q)}{\left(X,Y\right)}{II \left(X,Y\right)=\pi \left(\nabla 
^{M}_{X}Y\right).}
We define basic $r$-forms by :
$$\Omega_{B}^{r}\left(\mathcal{F}\right)=\left\{\Phi\in \Lambda ^{r}T^{*}M 
|\hskip0.2cm X\llcorner\Phi =0 \quad\text{and}\quad
X\llcorner d\Phi 
=0,\hskip0.5cm  \forall X \in \Gamma(L) \right\}, $$
where $d$ is the exterior derivative  and 
$X\llcorner$ is the interior product. We denote by  $d_{B}=d|_{\Omega _{B}\left(\mathcal{F}\right)}$
where $\Omega _{B}\left(\mathcal{F}\right)=\oplus_{r=0}^{p+q} \,\,\,\Omega_{B}^{r}\left(\mathcal{F}\right)$ and $\delta _{B}$  the adjoint operator of  $d_{B}$ with respect to the induced
scalar product. The basic Laplacian is defined as 
$ \Delta_{B}=d_{B}\delta_{B}+\delta_{B}d_{B}.$  Now we prove the following theorem.
 \begin{thm} \label{thm:06} Let $\left(M, g_{M},\mathcal{F}\right)$ be a compact Riemannian manifold with a Riemannian  foliation $\mathcal{F}$  and  a bundle-like metric $g_{M}$  with a coclosed basic $1$-form $\kappa.$ Assume that the transversal Ricci curvature is strictly positive, then the mean curvature $\kappa$ vanishes. 
\end{thm}
 {\bf Proof.} In \cite{h, MRT}, it is proved that the positivity of the transversal Ricci curvature implies the existence of a basic function $h$ such that $\kappa=d_{B}h.$ Then $\Delta_{B}h=\delta_{B}d_{B}h=\delta_{B}\kappa=0.$  Hence the harmonicity of $h$ implies that the function $h$ is closed, since $M$ is compact. Thus  the foliation is minimal. 
\hfill$\square$\\\\
Now, we assume that the normal bundle $Q$ carries a spin structure  and we denote  by $S(\mathcal{F})$  the  foliated spinor bundle.  The normal bundle acts on the spinor bundle by Clifford  multiplication  and the transversal Dirac operator \cite{gk1} is locally given by:
\begin{equation}
D_{tr}\Psi =\sum_{i=1}^{q}e_{i}\cdot\nabla  _{e_{i}}\Psi -\frac{1}{2} 
\kappa \cdot\Psi,
\label{eq:369}
\end{equation}
 for all $\Psi \in \Gamma(S(\mathcal{F})).$ We can easily prove using Green's theorem \cite{yt} that this operator is formally self-adjoint. We define the subspace of basic sections $\Gamma_{B}(S(\mathcal{F}))$ by 
 \begin{equation*}
 \Gamma_{B}(S(\mathcal{F}))= \{\Psi \in \Gamma(S(\mathcal{F})) |\, \nabla _{X} \Psi=0, \hskip0.5cm \forall X \in \Gamma(L)\}. 
 \end{equation*}
 The transversal Dirac operator leaves $\Gamma_{B}(S(\mathcal{F}))$ invariant if and only if the foliation is isoparametric. Moreover the basic Dirac operator defined by $D_b=D_{tr}|\Gamma_{B}(S(\mathcal{F})),$  has a discrete spectrum \cite{A} and  if the foliation $\mathcal{F}$ 
is isoparametric with $\delta_{B} 
\kappa =0$,   we have  the Schr{\"o}dinger-Lichnerowicz formula for $D_{b}$ \cite{gk1}
$$D_{b}^{2}\Psi =\nabla^{\star }\nabla\Psi 
+\frac{1}{4}K_{\sigma}^{\nabla}\Psi, $$
where $ K_{\sigma  }^{\nabla   }=\sigma^{\nabla}+|\kappa|^{2} $ and
$$ \nabla^{\star }\nabla\Psi =-\sum_{i=1}^{q}\nabla _{e_{i},e_{i}}^{2}\Psi 
+\nabla _{\kappa}\Psi, $$ with $
\nabla_{X,Y}^{2}=\nabla_{X}\nabla_{Y}-\nabla_{\nabla_{X}Y},$ for all $X,Y
\in \Gamma(TM).$ 
\section{Quaternion-K\"ahler Foliations}
\setcounter{equation}{0}
In this section,  we  review some basic relations on quaternion-K\"ahler spin foliations \cite{HM3} also we give  basic ingredients  for the estimate which  could be found in \cite{HM4} .\\
A foliation $\mathcal{F}$ of codimension $q=4m$  is said to be quaternion-K\"ahler if its principal bundle of oriented orthonormal frames $\SO Q$ admits a reduction $P$ to the subgroup $\Sp_{1} \cdot \Sp_{m}:=\Sp_{1} \times_{\Z_{2}} \Sp_{m} \subset \SO_{4m}.$ This is equivalent to the existence of  a subbundle $E$ of End($Q$)  of rank $3$ which admits  a local frame  $\{J_{\alpha}\}_{\alpha=1,2,3}$ such that the metric $g_{Q}$ is hermitian  for $J_{\alpha}, \alpha=1,2,3$ and verifies  
 \begin{equation}
\left\{\begin{array}{llc}
 J_{\alpha}\circ J_{\beta}=-\delta_{\alpha \beta} {\rm{Id}} + \varepsilon_{\alpha\beta\gamma}^{123}J_{\gamma}, \\
\nabla J_{\alpha}=\sum_{\beta=1}^{3}\omega_{\alpha}^{\beta}J_{\beta,}
\end{array}\right.
\label{eq:259}
\end{equation}
where $\omega_{\alpha}^{\beta}$ are the local $1$-forms  on $M$ and $\varepsilon_{\alpha \beta\gamma}^{123}=\pm 1$ if
$(\alpha,\beta,\gamma)$ is even or odd permutation of $(1,2,3).$  We note that a quaternion-K{\"a}hler foliation is  transversally Einstein \cite{Bes}, hence it admits a constant scalar curvature which is supposed to be positive throughout this paper.  
 A consequence of the definition is the existence of a parallel  4-form  $\Omega$ defined by 
$ \Omega =\sum_{\alpha=1}^{3} \Omega_{\alpha} \wedge \Omega_{\alpha},$
where the  $\Omega_{\alpha}$ are the local K\"ahler 2-forms associated with $J_{\alpha}.$
 The $4$-form  $\Omega$ can be written as
 \begin{equation}
 \Omega=\sum_{\alpha=1}^{3}\Omega_{\alpha}\cdot\Omega_{\alpha}+6m \rm Id.
\label{eq:558}
\end{equation}
Under the action of $\Omega,$  the foliated spinor bundle $S(\mathcal{F})$ splits into an orthogonal sum  $$ S(\mathcal{F})=\mathop\oplus_{r=0}^{m}S_{r}(\mathcal{F}),$$ where
$S_{r}(\mathcal{F})$ is the eigenbundle associated with the eigenvalue $\break \mu_{r}=6m-4r(r+2)$ of $\Omega.$ Moreover,  the action of the group $\Sp_{1}\times\Sp_{m}$ splits the bundle $Q^{\C}\otimes S_{r}(\mathcal{F})$ into \cite{ksw3} 
 \begin{eqnarray}
 Q^{\C}\otimes S_{r}(\mathcal{F})&=&W_{r+1,\bar{r}}(\mathcal{F})\oplus W_{r-1,\bar{r}}(\mathcal{F})\oplus W_{r+1,r-1}(\mathcal{F})\oplus W_{r-1,r+1}(\mathcal{F}) \nonumber \\ 
 && \oplus W_{r-1,r-1}(\mathcal{F})\oplus W_{r+1,r+1}(\mathcal{F}),
 \label{eq:1}
 \end{eqnarray} 
 where $W_{r,s}(\mathcal{F})$ denotes the space of the irreducible representation of the group $\Sp_{1}\times\Sp_{m}$ with dominant weight $$(r,1,\cdots,1,\underbrace {0,\cdots,0}_{s}), $$ and $W_{r,\bar{s}}(\mathcal{F})$ is the space of the irreducible representation of the group $\break \Sp_{1}\times\Sp_{m}$ with dominant weight $$(r,2,1,\cdots,1,\underbrace {0,\cdots,0}_{s}).$$ The last two bundles in \eqref{eq:1} are respectively isomorphic to $S_{r-1}(\mathcal{F})$ and $S_{r+1}(\mathcal{F}).$ We denote by $m_{r}$
the  restriction of the Clifford multiplication to  $\break Q^{\C}\otimes S_{r}(\mathcal{F}).$ The kernel of  $m_{r}$  splits into an orthogonal sum 
 \begin{equation*}
 \kker m_{r} = W_{r+1,\bar{r}}(\mathcal{F})\oplus W_{r-1,\bar{r}}(\mathcal{F})\oplus W_{r+1,r-1}(\mathcal{F})\oplus W_{r-1,r+1}(\mathcal{F}).
 \end{equation*}
 This comes from the computation of the image of $m_{r}$ of the maximal vector of each component of \eqref{eq:1}. Thus the restriction of $m_{r}$ to $ W_{r-1,r-1}(\mathcal{F})$ (resp. $W_{r+1,r+1}(\mathcal{F})$) is an isomorphism onto $S_{r-1}(\mathcal{F})$ (resp. $ S_{r+1}(\mathcal{F})$). Let (.\,,\,.)  be the usual hermitian product on $Q^{\C}\otimes S(\mathcal{F}).$ Since  $(m_{r}(.),m_{r}(.))$ and (.\,,\,.) are ($\Sp_{1}\times\Sp_{m}$)-invariant scalar products on  both $W_{r-1,r-1}(\mathcal{F})$ and $W_{r+1,r+1}(\mathcal{F}),$ one gets from Schur lemma 
\begin{equation} 
\forall w \in W_{r-1,r-1}(\mathcal{F}), \hskip0.2cm |m_{r}(w)|^{2}= \frac{2(r+1)(m-r+1)}{r}|w|^{2},
\label{eq:2}
\end{equation}
and,  
\begin{equation}
\forall w \in W_{r+1,r+1}(\mathcal{F}), \hskip0.2cm |m_{r}(w)|^{2}= \frac{2(r+1)(m+r+3)}{r+2}|w|^{2}.
\label{eq:3}
\end{equation} 
 In order to obtain a similar result for the other terms in \eqref{eq:1},  we locally  define   the operator  $\widetilde m : \Gamma(Q^{\C}\otimes S(\mathcal{F})) \longrightarrow \Gamma(E^{\C}\otimes S(\mathcal{F}))$  by  
 \begin{equation}
 \widetilde m(X\otimes\Psi)=\sum_{\alpha=1}^{3}J_{\alpha}\otimes(J_{\alpha}(X)\cdot \Psi), 
 \label{eq:111}
 \end{equation}
 for all $X\in \Gamma(Q)$ and $\Psi\in\Gamma(S(\mathcal{F})).$ We  denote by $\widetilde m_{r}$  the restriction of $\widetilde m$ to $Q^{\C}\otimes S_{r}(\mathcal{F}).$ As above, computing the image of  $\widetilde m_{r}$ of maximal vector of each component of \eqref{eq:1}, the kernel of $\widetilde m_{r}$  splits into
\begin{equation*}
\kker \widetilde m_{r}=W_{r+1,\bar{r}}(\mathcal{F})\oplus W_{r-1,\bar{r}}(\mathcal{F}).
\end{equation*}
Using the same argument as in \eqref{eq:2} and \eqref{eq:3}, one gets from Schur lemma
\begin{eqnarray}
\forall w  \in W_{r+1,r-1}(\mathcal{F}), \hskip0.2cm |\widetilde m_{r}(w)|^{2} &=& 4(m-r+1) |w|^{2},\label{eq:6}\\
\forall w  \in W_{r-1,r+1}(\mathcal{F}),\hskip0.2cm |\widetilde m_{r}(w)|^{2} &=& 4(m+r+3) |w|^{2},\label{eq:7}\\
\forall w  \in W_{r-1,r-1}(\mathcal{F}), \hskip0.2cm |\widetilde m_{r}(w)|^{2} &=& \frac{2(r-1)(m-r+1)}{r} |w|^{2},\label{eq:8}\\
\forall w  \in W_{r+1,r+1}(\mathcal{F}), \hskip0.2cm  |\widetilde m_{r}(w)|^{2} &=& \frac{2(r+3)(m+r+3)}{r+2} |w|^{2}.\label{eq:9}
\end{eqnarray}
\section{The main Result}
\setcounter{equation}{0}
In this section, we show \eqref{eq:16} by using the decomposition of the bundle $Q^{\C}\otimes S_{r}(\mathcal{F})$ given in the above section . We refer to \cite{HM4}, \cite{ksw1}, \cite{ksw3}.
\begin{thm} Under the same conditions as in  Theorem \emph{\ref{thm:06}} with the assumption that the foliation $\mathcal{F}$ has a quaternion-K\"ahler spin structure of  codimension $q=4m,$  then  the mean curvature $\kappa$ vanishes and any eigenvalue $\lambda$ of the basic Dirac operator satisfies 
\begin{equation*}
\lambda^{2}\geq\frac{m+3}{4(m+2)}\sigma^{\nabla}, 
\end{equation*}
where $\sigma^{\nabla}$ denotes the transversal scalar curvature.
\end{thm}
{\bf Proof.} The fact that $\mathcal{F}$ is minimal comes from Theorem \ref{thm:06} since the transversal scalar curvature is supposed to be positive.  For the second part, according to the decomposition  \eqref{eq:1},  for any  $\Psi \in \Gamma_B(S_{r}(\mathcal{F})),$ the covariant derivative $ \nabla\Psi$ splits into
\begin{equation}
\begin{array}{lcc}
\nabla\Psi=& (\nabla\Psi)_{r+1,\bar{r}}+(\nabla\Psi)_{r-1,\bar{r}}+(\nabla\Psi)_{r+1,r-1}+(\nabla\Psi)_{r-1,r+1}
\\\\&+ (\nabla\Psi)_{r-1,r-1}+(\nabla\Psi)_{r+1,r+1}.
\end{array}
\label{eq:69} 
\end{equation}
 In order to compute the norm of $\nabla\Psi,$ since the last two terms in the above equation are  sections in the subbundles $S_{r-1}(\mathcal{F})$ and  $S_{r+1}(\mathcal{F})$ respectively, we get from \eqref{eq:2} and \eqref{eq:3},
\begin{eqnarray} 
|(\nabla\Psi)_{r-1,r-1}|^{2}=\frac{r}{2(r+1)(m-r+1)}|D_{-}\Psi|^{2},
\label{eq:22}
\end{eqnarray}
and, 
\begin{eqnarray}
|(\nabla\Psi)_{r+1,r+1}|^{2}=\frac{r+2}{2(r+1)(m+r+3)}|D_{+}\Psi|^{2},
\label{eq:23}
\end{eqnarray}
where  $D_{-}\Psi=(D_{b}\Psi)_{r-1}$ and   $D_{+}\Psi=(D_{b}\Psi)_{r+1}.$ Similar results could be obtained for the other terms in \eqref{eq:69} by using the definition of the operator $\widetilde m$  in \eqref{eq:111}. For this, we consider for any spinor $\Psi$ the operator $D_{\alpha}\Psi$  locally defined by  $\sum_{i=1}^{4m}J_{\alpha}(e_{i})\cdot\nabla_{e_{i}}\Psi.$  Hence we have $\widetilde m(\nabla\Psi)=\sum_{\alpha=1}^{3}J_{\alpha}\otimes D_{\alpha}\Psi$ and we get that 
\begin{equation*}
|\widetilde{m}(\nabla\Psi)|^{2}=\sum_{\alpha=1}^{3}|D_{\alpha}\Psi|^{2}.
 \end{equation*}
 On the other hand, 
Equations \eqref{eq:6}, \eqref{eq:7}, \eqref{eq:8}, \eqref{eq:9} imply that
\begin{eqnarray*}
|\widetilde m((\nabla\Psi)_{r+1,r-1})|^{2} &=& 4(m-r+1)|(\nabla\Psi)_{r+1,r-1}|^{2},\\
|\widetilde m((\nabla\Psi)_{r-1,r+1})|^{2} &=& 4(m+r+3)|(\nabla\Psi)_{r-1,r+1}|^{2},\\
|\widetilde m((\nabla\Psi)_{r-1,r-1})|^{2} &=& \frac{2(r-1)(m-r+1)}{r}|(\nabla\Psi)_{r-1,r-1}|^{2},\\
|\widetilde m((\nabla\Psi)_{r+1,r+1})|^{2} &=& \frac{2(r+3)(m+r+3)}{r+2}|(\nabla\Psi)_{r+1,r+1}|^{2}.
\end{eqnarray*}
Hence by the above equations and  \eqref{eq:22}, \eqref{eq:23}, we conclude  for any  $\Psi \in \Gamma_B(S_{r}(\mathcal{F}))$ that
\begin{eqnarray}
\sum_{\alpha=1}^{3}|D_{\alpha}\Psi|^{2}&=& 4(m-r+1)|(\nabla\Psi)_{r+1,r-1}|^{2}+4(m+r+3)|(\nabla\Psi)_{r-1,r+1}|^{2}\nonumber\\&&+\frac{r+3}{r+2}|D_{+}\Psi|^{2}+\frac{r-1}{r+1}|D_{-}\Psi|^{2}.
\label{eq:55}
\end{eqnarray}
Then using equations \eqref{eq:22}, \eqref{eq:23}, \eqref{eq:55} and by \eqref{eq:69},  we write the norm of  $\nabla\Psi$  as 
\begin{eqnarray}
|\nabla\Psi|^{2} &=& |(\nabla\Psi)_{r+1,\bar {r}}|^{2}+|(\nabla\Psi)_{r-1,\bar{r}}|^{2} 
+\frac{2(r+1)}{m+r+3}|(\nabla\Psi)_{r+1,r-1}|^{2} \nonumber\\&& +\frac{1}{4(m+r+3)}\sum_{\alpha=1}^{3}|D_{\alpha}\Psi|^{2}+\frac{1}{4(m+r+3)}|D_{+}\Psi|^{2} \nonumber\\
&&+\frac{m+3r+1}{4(m-r+1)(m+r+3)}|D_{-}\Psi|^{2}.
\label{eq:70}
\end{eqnarray}
Now  let $\lambda$ be any eigenvalue of the basic Dirac operator, then there exists an eigenspinor $\Psi,$ called of type ($r,r+1$), such that
\begin{equation*}
D_{b}\Psi=\lambda \Psi  \quad\text{and}\quad  \Psi= \Psi_{r}+\Psi_{r+1},
\end{equation*}
 with $r\in \{0,\cdots,m-1\}.$ In \cite{HM3}, it is showed that for any spinor $\break \Psi \in \Gamma_B(S(\mathcal{F})),$ we have 
\begin{equation*}
\int_{M}\sum_{\alpha=1}^{3}|D_{\alpha}\Psi|^{2}=3\int_{M}(D_b^{2}\Psi,\Psi)+\frac{\sigma^\nabla}{4m(m+2)}\int_{M}((\Omega-6m)\cdot\Psi,\Psi).
\end{equation*}
Therefore, applying Equation \eqref{eq:70} to $\Psi_{r+1}$ and integrating  over $M,$ one gets since $D_{-}\Psi_{r+1}=\lambda\Psi_{r}$ and $D_{+}\Psi_{r+1}=0$
\begin{eqnarray*}
0 &\leq& ||\nabla\Psi_{r+1}||_{L^{2}} ^{2}-a_{r}\lambda^{2}||\Psi_{r+1}||_{L^{2}}^{2}+b_{r} \sigma^{\nabla}||\Psi_{r+1}||_{L^{2}}^{2} -c_{r}\lambda^{2}||\Psi_{r}||_{L^{2}}^{2},
\end{eqnarray*}
 where, 
 \begin{equation*}
 \left\{\begin{array}{llc}
 a_{r}=\frac{3}{4(m+r+4)}, \\\\
 b_{r}=\frac{(r+1)(r+3)}{4m(m+2)(m+r+4)},\\\\
 c_{r}=\frac{m+3r+4}{4(m-r)(m+r+4)}. 
 \end{array}\right.
 \end{equation*}
 Finally with the help of  the Schr\"odinger-Lichnerowicz formula  and  the fact that $\Psi_{r}$ and $\Psi_{r+1}$ have the same  $L^{2}$-norms, we get \eqref{eq:16}.
\hfill$\square$
\section{ The Limiting  case} \label{sec:4}
 \setcounter{equation}{0}
 Let $\lambda$ be the first eigenvalue satisfying equality in \eqref{eq:16} and $\Psi$ an eigenspinor of type $(r,r+1).$  From the proof of  Theorem \ref{thm:06}, one gets necessarily that $r=0$ and  the following equations \cite{HM4} 
 \begin{equation}
 \left\{\begin{array}{llc}
 |\nabla\Psi_{0}|^2=\frac{1}{m+3}|D_{b}\Psi_{0}|^{2},\\\\
 |\nabla\Psi_{1}|^2=\frac{1}{4m}|D_{b}\Psi_{1}|^{2}+\frac{1}{4(m+4)}\sum_{\alpha=1}^{3}|D_{\alpha}\Psi_{1}|^{2}.
 \end{array}\right.
 \label{eq:599}
 \end{equation}
 Furthermore,  the spinor  $\Psi_{1}$ satisfies 
 \begin{equation}
 \begin{array}{llc}
 \sum_{\alpha=1}^{3}\Omega_{\alpha}\cdot D_{\alpha}\Psi_{1}=0, \\\\
 \sum_{\beta, \gamma} \varepsilon_{\alpha \beta \gamma}^{123} \Omega_{\beta} \cdot D_{\gamma}\Psi_{1}=8D_{\alpha}\Psi_{1}, & \textrm{$\forall \alpha=1,2,3.$}
 \label{eq:889}
 \end{array}
 \end{equation}
 Moreover for all $X\in \Gamma(Q),$  we have  the {\it quaternion-K\"ahler Killing} equations \cite{HM4}, \cite{ksw1}, \cite{ksw2}
 \begin{equation}
 \nabla_{X}\Psi_{0}=-\frac{\lambda}{m+3}p_{1}(X)\cdot\Psi_{1}, 
 \label{eq:24}
 \end{equation}
 and, 
 \begin{equation}
 \nabla_{X}\Psi_{1}=-\frac{\lambda}{4m}X\cdot\Psi_{0} -\frac{1}{4(m+4)}\sum_{\alpha=1}^{3}J_{\alpha}(X)\cdot D_{\alpha}\Psi_{1}, 
 \label{eq:25}
 \end{equation}
 where for all $X\in \Gamma(Q),$  the operator $p_{1}$ is defined by (see \cite{HM1}) 
 \begin{equation*}
 \left\{\begin{array}{llc}
 p_{1}(X)=\frac{1}{8}(5X+\mathcal{J}(X)),\\\\  
  \mathcal{J}(X)=\frac{1}{4}[\Omega,X].
  \end{array}\right.
  \end{equation*}   
 In order to  prove \eqref{eq:24}, we  define the transversal quaternion-K\"ahler twistor operator,  denoted by $\mathcal{P}^0,$  on the bundle  $S_{0}(\mathcal{F})$ whose the image lies in the bundle $Q^{*}\otimes S_{0}(\mathcal{F})$ (see \cite{HM1} for the details). For any spinor field $\psi_0 \in \Gamma_B(S_{0}(\mathcal{F})),$  we write
 \begin{equation*}
 \mathcal{P}^0\psi_0=\sum_{i=1}^{4m} e_i\otimes(\nabla_{e_{i}}\psi_0+\frac{1}{m+3}p_{1}(e_{i})\cdot D_{b}\psi_0),
 \end{equation*}
  where $\{e_{i}\}_{i=1,\cdots,4m}$ is a local orthonormal frame of  $\Gamma(Q).$ By a straightforward computation and  with the definition of  $p_{1},$  we easily verify that $\break\sum_{i=1}^{4m} e_i\cdot \mathcal{P}^0_{e_{i}}\psi_0=0.$ Hence the image of  $\mathcal{P}^0$ lies in the kernel of Clifford multiplication $m_{0}$.  Since $ \mathcal{P}^0_{e_{i}}\psi_0$ is a section on $S_{0}(\mathcal{F}),$  we deduce with the definition of the operator $\mathcal{J},$ that  $\sum_{i=1}^{4m}\mathcal{J}(e_{i})\cdot \mathcal{P}^0_{e_{i}}\psi_0=0.$  Then 
 \begin{eqnarray}
 |\mathcal{P}^0\psi_{0}|^2&=&\sum_{i=1}^{4m}( \mathcal{P}^0_{e_{i}}\psi_0,\mathcal{P}^0_{e_{i}}\psi_0)\nonumber\\
 &=&\sum_{i=1}^{4m}( \mathcal{P}^0_{e_{i}}\psi_0, \nabla_{e_{i}}\psi_{0})\nonumber\\&=&|\nabla\psi_{0}|^2+\frac{1}{m+3} \sum_{i=1}^{4m}(p_{1}(e_{i})\cdot D_{b}\psi_{0},\nabla_{e_{i}}\psi_{0}).
 \label{eq:355}
 \end{eqnarray}
 Since Clifford multiplication by $\mathcal{J}$ is symmetric, one can easily verify that $(\nabla_{e_{i}}\psi_{0},p_{1}(e_{i}) \cdot D_{b}\psi_{0})=-(e_{i}\cdot\nabla_{e_{i}}\psi_{0},D_{b}\psi_{0}).$  Then for any spinor  $\break\psi_{0} \in \Gamma(S_{0}(\mathcal{F})),$  Equation \eqref{eq:355} reduces to 
 \begin{eqnarray*}
  |\mathcal{P}^0\psi_{0}|^2=|\nabla\psi_{0}|^2-\frac{1}{m+3}|D_{b}\psi_{0}|^{2},
 \end{eqnarray*}
which  vanishes by  \eqref{eq:599}  for the spinor field $\Psi_{0}.$  Thus Equation \eqref{eq:24} is satisfied for $X=e_{i}.$
\hfill$\square$\\\\
Now, we will prove Equation \eqref{eq:25}.  The proof consists in computing  the  sum
\begin{equation}
\begin{array}{llc}
\sum_{i=1}^{4m}|\nabla_{e_i}\Psi_1+\frac{1}{4m}e_{i}\cdot D_b \Psi_1+\frac{1}{4(m+4)}\sum_{\alpha=1}^{3}J_{\alpha}e_{i}\cdot D_{\alpha}\Psi_{1}|^2=\\\\
|\nabla\Psi_1|^2+\frac{1}{4m}|D_b \Psi_1|^2+\frac{1}{16(m+4)^2}\sum_{i=1}^{4m}|\sum_{\alpha=1}^{3}J_{\alpha}e_{i}\cdot D_{\alpha}\Psi_{1}|^2\\\\+\frac{1}{2m}\sum_{i=1}^{4m}(\nabla_{e_i}\Psi_1,e_{i}\cdot D_b \Psi_1)+\frac{1}{2(m+4)}\sum_{i=1}^{4m}(\nabla_{e_i}\Psi_1,J_{\alpha}e_{i}\cdot D_{\alpha}\Psi_{1}).
\label{eq:358}
\end{array}
\end{equation}
 The fact that Clifford multiplication by $e_{i}$ and $J_{\alpha}(e_{i})$ is skew-symmetric,  the last terms are easily computed and  it remains to compute the third  term in the r.h.s. of  \eqref{eq:358} . For this, using a local orthonormal frame $\{ J_{\alpha} e_{i} \}_{i=1,\cdots,4m}$ and  \eqref{eq:259}, it follows
\begin{eqnarray}
\sum_{i=1}^{4m}|\sum_{\alpha=1}^{3}J_{\alpha}e_{i}\cdot D_{\alpha}\Psi_{1}|^2&=&\sum_{i,\alpha,\beta} (J_{\alpha}e_{i}\cdot D_{\alpha}\Psi_{1},J_{\beta}e_{i}\cdot D_{\beta}\Psi_{1})\nonumber\\
&=&\sum_{i,\alpha,\beta}(D_{\alpha}\Psi_{1},e_{i} \cdot J_{\beta}J_{\alpha}e_{i}\cdot D_{\beta}\Psi_{1})\nonumber\\
&=&4(m+4)\sum_{\alpha=1}^3 |D_{\alpha}\Psi_{1}|^2.
\label{eq:20}
\end{eqnarray}
The last identity in \eqref{eq:20} comes from \eqref{eq:259} and \eqref{eq:889}.  Finally substituting \eqref{eq:20} and using \eqref{eq:369}, Equation \eqref{eq:358} reduces to 
$$
\begin{array}{llc}
\sum_{i=1}^{4m}|\nabla_{e_i}\Psi_1+\frac{1}{4m}e_{i}\cdot D_b \Psi_1+\frac{1}{4(m+4)}\sum_{\alpha=1}^{3}J_{\alpha}e_{i}\cdot D_{\alpha}\Psi_{1}|^2=\\\\
|\nabla\Psi_1|^2-\frac{1}{4m}|D_b \Psi_1|^2-\frac{1}{4(m+4)}\sum_{\alpha=1}^{3} |D_{\alpha}\Psi_{1}|^2,
\end{array}
$$
which  vanishes by \eqref{eq:599}. 
\hfill$\square$\\\\
{\bf Example 1} We consider the compact manifold $N=M\times \HH {\rm P}^{m},$  where $M$ is a compact Riemannian manifold of dimension $p$ and  $\HH {\rm P}^{m}$ is the quaternionic projective space with its standard metric. Let $g_{N}$ be the product metric  on $N.$  We define a foliation $\mathcal{F}$ on $N$ by its leaves of the form $M\times\{y\}$ where $y\in  \HH {\rm P}^{m}.$ This is a Riemannian foliation on $N$ and $g_{N}$ is a bundle-like metric with totally geodesic fibers. Since the fibers of the normal bundle are the tangent space of $\HH {\rm P}^{m},$ then it carries a quaternion-K\"ahler spin structure and  the basic Dirac operator coincides with the one on  $\HH {\rm P}^{m}$ where the eigenvalues are computed in \cite{JL}. Hence the limiting case in \eqref{eq:16} is achieved.\\\\
{\bf Example 2} Let $M$ be a  compact  $3$-Sasakian manifold and consider the foliation  on $M$ defined by its Killing vector fields.  This  is a Riemannian foliation with a bundle-like metric and totally geodesic fibers diffeomorphic to $\Gamma \setminus S^{3}$ where  $\Gamma$ is a finite subgroup of $\Sp_{1}$ \cite{BGM}.  It induces a quaternion-K\"ahler  spin structure on the normal bundle with positive transversal scalar curvature.  If  $M$ is either $S^{4q+3}$ or $\R \PP ^{4q+3},$ then it   projects onto  $\HH {\rm P}^{m}$ (Hopf fibration). Since the  fibers of the normal bundle are isomorphic to the tangent space of  $\HH {\rm P}^{m},$  then  equality in \eqref{eq:16} is achieved.

\end{document}